\documentclass[10pt]{article}

\usepackage{amscd,amsmath, amssymb, fancyhdr, graphics, url}

\usepackage[backref=page]{hyperref}
\renewcommand*{\backref}[1]{}
\renewcommand*{\backrefalt}[4]{%
    \ifcase #1 (Not cited.)%
    \or        (Cited on page~#2.)%
    \else      (Cited on pages~#2.)%
    \fi}

\hypersetup{
     colorlinks   = true,
     citecolor    = magenta,
     linkcolor= blue
}

\numberwithin{equation}{section}

\newcommand{\version}{version 2.1,\ \ September 22, 2019}

\makeatletter
\def\x@arrow{\DOTSB\Relbar}
\def\xlongrightarrowfill@{\arrowfill@\relbar\relbar\longrightarrow}
\newcommand{\xlongrightarrow}[2][]{%
        \ext@arrow 0099\xlongrightarrowfill@{#1}{#2}}
\makeatother

\def\eqref#1{(\ref{#1})}
\newcommand{\goth}{\mathfrak}

\newcommand{\arrow}{{\:\longrightarrow\:}}
\newcommand{\Z}{{\Bbb Z}}
\def\C{{\Bbb C}}
\newcommand{\R}{{\Bbb R}}

\def\1{\sqrt{-1}\:}

\newcommand{\cntrct}                
{\hspace{2pt}\raisebox{1pt}{\text{$\lrcorner$}}\hspace{2pt}}

\renewcommand{\tilde}{\widetilde}

\renewcommand{\phi}{\varphi}
\renewcommand{\epsilon}{\varepsilon}
\renewcommand{\geq}{\geqslant}

\newcommand{\Gr}{\operatorname{Gr}}
\newcommand{\Teich}{\operatorname{\sf Teich}}
\newcommand{\Comp}{\operatorname{\sf Comp}}

\newcommand{\Per}{\operatorname{\sf Per}}
\newcommand{\Perspace}{\operatorname{{\Bbb P}\sf er}}

\newcommand{\Vol}{\operatorname{Vol}}

\newcommand{\Aut}{\operatorname{Aut}}
\newcommand{\Mon}{\operatorname{\sf Mon}}

\newcommand{\Diff}{\operatorname{\sf Diff}}

\renewcommand{\Re}{\operatorname{Re}}
\renewcommand{\Im}{\operatorname{Im}}

\newcounter{Mycounter}[section]
\newcounter{lemma}[section]
\renewcommand{\thelemma}{{Lemma \thesection.\arabic{lemma}}}
\newcommand{\lemma}{%
    \setcounter{lemma}{\value{Mycounter}}
    \refstepcounter{lemma}
    \stepcounter{Mycounter}
    {\noindent \bf \thelemma:\ }}

\newcounter{claim}[section]
\renewcommand{\theclaim}{{Claim \thesection.\arabic{claim}}}
\newcommand{\claim}{%
    \setcounter{claim}{\value{Mycounter}}
    \refstepcounter{claim}
    \stepcounter{Mycounter}
    {\noindent \bf \theclaim:\ }}

\newcounter{sublemma}[section]
\newcounter{corollary}[section]
\renewcommand{\thecorollary}{{Corollary \thesection.\arabic{corollary}}}
\newcommand{\corollary}{%
    \setcounter{corollary}{\value{Mycounter}}
    \refstepcounter{corollary}
    \stepcounter{Mycounter}
    {\noindent \bf \thecorollary:\ }}

\newcounter{theorem}[section]
\renewcommand{\thetheorem}{{Theorem \thesection.\arabic{theorem}}}
\newcommand{\theorem}{%
    \setcounter{theorem}{\value{Mycounter}}
    \refstepcounter{theorem}
    \stepcounter{Mycounter}
    {\noindent \bf \thetheorem:\ }}

\newcounter{conjecture}[section]
\newcounter{proposition}[section]
\newcounter{definition}[section]
\renewcommand{\thedefinition}
      {{Definition~\thesection.\arabic{definition}}}
\newcommand{\definition}{%
    \setcounter{definition}{\value{Mycounter}}
    \refstepcounter{definition}
    \stepcounter{Mycounter}
    {\noindent \bf \thedefinition:\ }}

\newcounter{example}[section]
\newcounter{remark}[section]
\renewcommand{\theremark}{{Remark \thesection.\arabic{remark}}}
\newcommand{\remark}{%
    \setcounter{remark}{\value{Mycounter}}
    \refstepcounter{remark}
    \stepcounter{Mycounter}
    {\noindent \bf \theremark:\ }}

\newcounter{problem}[section]
\newcounter{question}[section]
\newcommand{\proof}{{\bf Proof:\:}}

\makeatletter

\@addtoreset{equation}{section} \@addtoreset{footnote}{section}
\makeatother

\def\blacksquare{\hbox{\vrule width 5pt height 5pt depth 0pt}}
\def\endproof{\blacksquare}

\begin{document}

\begin{center}
{\LARGE\bf
Mapping class group\\[1mm] and a global Torelli theorem\\[1mm] for hyperk\"ahler manifolds:\\[1mm] an erratum\\[4mm]
}

 Misha
Verbitsky\footnote{Partially supported by
by  the  Russian Academic Excellence Project '5-100', FAPERJ E-26/202.912/2018 
and CNPq - Process 313608/2017-2


}

\end{center}

{\small \hspace{0.05\linewidth}
\begin{minipage}[t]{0.95\linewidth}
{\bf Abstract} \\
{\bf A mapping class group} of an oriented manifold
is a quotient of its diffeomorphism group by the  isotopies.
In the published version of ``Mapping class group and a
global Torelli theorem for hyperk\"ahler
manifolds'' I made an error based on a
wrong quotation of Dennis Sullivan's famous paper
``Infinitesimal computations in topology''. I claimed
that the natural homomorphism from the mapping class
group to the group of automorphims of cohomology
of a simply connected K\"ahler manifold has finite kernel.
In a recent preprint \cite{_Kreck_Su_}, Matthias Kreck and Yang Su produced
counterexamples to this statement.
Here I correct this error and other related errors,
observing that the results of
``Mapping class group and a global Torelli theorem''
remain true after an appropriate change of terminology.
\end{minipage}
}

\tableofcontents

\hfill

Throughout this paper, ``hyperk\"ahler manifold''
means compact hyperk\"ahler manifold of maximal holonomy, 
also known as IHS (irreducible holomorphically symplectic) manifold.

\section{Mapping class group}
\label{_MCG_Section_}

Let $M$ be a manifold, $\Diff(M)$ its diffeomorphism
group (sometimes - its oriented diffeomorphism group), 
and $\Diff_0(M)$ the connected component of unity in
$\Diff(M)$. The quotient group
$\Gamma=\frac{\Diff(M)}{\Diff_0(M)}$
is called {\bf the mapping class group}.
In \cite{_Sullivan:infinite_}, Dennis Sullivan
compared the mapping class group to the algebraic
group $A$ of automorphisms of the minimal model
of the de Rham differential graded algebra
preserving the integer structure and Pontryagin classes, and found
that $\Gamma$ is mapped to a lattice of integer points of 
$A$ with finite kernel and cokernel whenever $M$ has nilpotent
fundamental group and $\dim M \geq 5$.

As shown in \cite{_DGMS:Formality_}, 
the de Rham algebra of a K\"ahler manifold $M$
is formal, that is, quasi-isomorphic to its cohomology. 
I claimed that this implies that the group $A$ is isomorphic to the automorphisms 
of the cohomology preserving the Pontryagin classes.

This is wrong, because in the typical situation 
the minimal model is still much bigger than the cohomology algebra,
and its automorphism group is also bigger. 
Indeed, as shown by Kreck and Su 
(\cite{_Kreck_Su_}), the natural map from
the mapping class group to $\Aut(H^*(M))$
can have infinitely-dimensional kernel.

Therefore, Theorem 3.5 from \cite{_V:Torelli_} is false as
stated. The correct version is the following.

\hfill

\theorem \label{_mapping_class_intro_Equation_}
Let $M$ be a compact hyperk\"ahler manifold of maximal holonomy, and
$\Gamma= \Diff/\Diff_0$ its mapping class group.
Then $\Gamma$ acts on $H^2(M, \R)$ preserving the
Bogomolov-Beauville-Fujiki form. Moreover, the corresponding
homomorphism $\phi:\; \Gamma \arrow O(H^2(M, \Z), q)$  has finite index in
$O(H^2(M, \Z), q)$. 

\hfill

\proof
The image of  $\phi$ is described in \cite[Theorem 3.5]{_V:Torelli_},
using the results of hyperk\"ahler geometry (the local Torelli theorem,
Riemann-Hodge reltions and and the computation of cohomology algebra obtained in 
\cite{_Verbitsky:coho_announce_}). The kernel (claimed to be finite in
\cite[Theorem 3.5]{_V:Torelli_}) is not finite. 
\endproof

\hfill

The correct 
quotation from Sullivan's paper \cite[Theorem 13.3]{_Sullivan:infinite_} 
follows.

\hfill

\definition\label{_algebraic_MCG_Definition_}
Let $\cal M$ be the minimal model of de Rham algebra of a 
manifold $M$, $\sigma$ its automorphism, and $\xi\in \cal M$
a class which satisfies $d\xi = \sigma(p)-p$, where
$p=(p_1, p_2, ..., p_n)$ is the total Pontryagin class.
The pair $(\sigma, \xi)$ is called {\bf an algebraic
diffeomorphism} of $M$. An algebraic diffeomorphism
$(\sigma, \xi)$ is called {\bf an algebraic isotopy} if
there exists a derivation $\delta:\; {\cal M} \arrow {\cal M}$ 
of degree -1 such that $\sigma = e^{d\delta + \delta d}$ and
\[ \xi=\delta\left(\sum_{i=1}^\infty \frac{(d\delta)^i}{(i+1)!}(p)\right)
\]
is an explicit homotopy between $p$ and $\sigma(p)$ induced by
the expression $\sigma = e^{d\delta + \delta d}$.
There is a natural group structure on the set of algebraic
diffeomorphisms ${\cal D}$ and algebraic isotopies ${\cal I}$;
moreover, ${\cal I}$ is a normal subgroup of ${\cal D}$.
{\bf An algebraic mapping class group} is the 
quotient $\frac{\cal D}{\cal I}$. It is (generally speaking)
a pro-algebraic group, with the natural integer structure
induced from the integer structure on the cohomology of 
the de Rham algebra. 

\hfill

\theorem (\cite[Theorem 13.3]{_Sullivan:infinite_})
Let $M$ be a compact, simply-connected manifold of dimension $\geq 5$.
Then its mapping class group is commensurable with
the group of integer points in its algebraic mapping class group
$\frac{\cal D}{\cal I}$.

\hfill

\remark If the manifold
$M$ is, in addtion, K\"ahler, its de Rham algebra is formal.
Then the Sullivan's minimal model ${\cal M}$ can be constructed from its
cohomology algebra $H^*(M)$. Any automorphism of 
$H^*(M)$ gives an automorphism of ${\cal M}$, but
 ${\cal M}$ can {\em a priori} have more automorphisms
than $H^*(M)$. This explains the error in 
\cite[Theorem 3.5]{_V:Torelli_}.

\hfill

\remark \label{_Torelli_group_Remark_}
Define {\bf the Torelli group} ${\cal T}$
 of $M$ as the subgroup of all elements of the
mapping class group of $M$ acting trivially on $H^*(M)$.
Kreck and Su (\cite{_Kreck_Su_}) have computed the Torelli
group for two 8-dimensional examples of hyperk\"ahler manifolds: 
the generalized Kummer variety and the second Hilbert scheme of K3 surface.
The Torelli group is infinite for the first one, and finite for the second.

\hfill  

\remark Sometimes one defines the Torelli group
as a subgroup $\tilde {\cal T}$ of the mapping class group which acts trivially
on the cohomology. In \cite[Theorem 3.5 (ii)]{_V:Torelli_}, it was shown that
the homomorphism $\Aut(H^*(M,\Z), O(H^2(M,\Z))$ has finite
kernel, hence $\tilde {\cal T}$ has finite index in ${\cal T}$.


\section{The monodromy group and the marked moduli space}
\label{_marked_Section_}

Let $M$ be a manifold and $\Comp$ the set of all complex
structures on $M$. Interpreting a complex structure as an
endomorphism of the tangent space, we imbue $\Comp$
with $C^\infty$-topology. Let $\Diff_0$ be the isotopy
group (that is, connected component of the unity in the
diffeomorphism group), and $\Gamma:= \frac{\Diff}{\Diff_0}$
its mapping class group. {\bf The Teichm\"uller space of complex structures on $M$}
is the quotient $\Teich:=\frac{\Comp}{\Diff_0}$ equipped
with the quotient topology. In many cases
(for compact Riemann surfaces, Calabi-Yau manifolds,
compact tori and hyperk\"ahler manifolds) the 
Teichm\"uller space is smooth. Its quotient
$\Teich/\Gamma$
parametrizes the complex structures on $M$, and 
(in many situations) can be understood as the moduli
space of deformations. 

\hfill

For convenience, we redefine the notation for
hyperk\"ahler manifolds. In this case, $\Teich$ denotes
the Teichm\"uller space of all complex structures 
of hyperk\"ahler type (that is, complex structures which are 
holomorphically symplectic
and K\"ahler). From the local Torelli theorem it follows
that this space is open in the more general Teichm\"uller
space of all complex structures.

\hfill

\claim\label{_Torelli_action_Claim_}
Let $K\subset \Gamma$ be the Torelli group of $M$, that
is, the group of all elements trivially acting on $H^2(M)$.
Consider an element $\gamma \in K$  which fixes a point
$I\in \Teich$. Then $\gamma$ acts trivially 
on the corresponding connected component of $\Teich$.
Moreover, the group $\goth S_I$ of such $\gamma$ is finite.

\hfill

\proof
Let $I\in \Teich$ be a fixed point of an element $\gamma\in \Gamma$,
and $l\in \Gr_{+,+}(H^2(M, \R))$ be the 2-plane
in the Grassmann space of positive 2-planes in $H^2(M,\R)$ corresponding to $I$, with 
$l=\langle \Re\Omega, \Im \Omega)$, where $\Omega\in H^2(M, \C)$
is the cohomology class of the holomorphic symplectic
form. Each connected component of the Hausdorff reduction
of $\Teich$ is naturally identified with  the space
$\Gr_{+,+}(H^2(M,\R))$ by the global Torelli theorem.

Let $\gamma$ be an element of the Torelli group of
$M$ which fixes $I\in\Teich$. Since $\gamma$
commutes with the period map $\Per:\;\Teich^I \arrow \Gr_{+,+}(H^2(M,\R))$
and acts trivially on $\Gr_{+,+}(H^2(M,\R))$, it maps
the set $\Per^{-1}(l)\subset \Teich^I$ to itself, for each $l\in \Perspace$.
However, $\Per^{-1}(l)$ is identified with the set
of all K\"ahler chambers in $H^{1,1}(M,I)$ for a complex structure
$I\in \Per^{-1}(l)$ (\cite[Theorem 5.16]{_Markman:survey_}). 
Since $\gamma$ acts trivially on $H^2(M)$,
it acts trivially on the set of K\"ahler chambers.
Then it acts trivially on the whole $\Teich^I$.

It remains to show only that $\goth S_I$ is finite.
Let $\Diff\cdot I$ be the orbit of $I$,
and $G$ the stabiliser of $I$, that is,
the group of complex automorphisms of $(M,I)$.
Then $\Diff\cdot I= \Diff/G$. Therefore, the set
of connected components of the orbit $\Diff\cdot I$
is parametrized by $\Gamma/G$. We obtain that 
$\gamma\in \Gamma$ fixes a point $J$ in $\Teich$ if
and only if it acts on $(M,J)$ by complex automorphisms.

Using Calabi-Yau theorem, it is easy to prove that 
the group of complex automorphisms acting trivially on
$H^2(M)$ is finite (\cite[Theorem 4.26]{_V:Torelli_}). 
Indeed, any such automorphism
preserves the Calabi-Yau metric, which is uniquely
determined by its K\"ahler class, and the group
of isometries of a compact metric space is compact.
\endproof

\hfill

\definition
Let $M$ be a hyperk\"ahler manifold.
The {\bf marked moduli space} of $M$ is 
the quotient $\Teich/K$.
From \ref{_Torelli_action_Claim_} it follows that
each connected component of $\Teich/K$ is diffeomorphic
to the corresponding connected component of $\Teich$,
because $K$ acts by permuting isomorphic connected components
of $\Teich$.

\hfill

\definition
Let  $\Teich^I$ denote the connected component of
$\Teich$ containing a point $I\in \Teich$. 
{\bf Monodromy group} $\Mon_I$ of a hyperk\"ahler manifold
$(M,I)$ is the group of all mapping class elements
preserving the component $\Teich^I$.

\hfill

\remark
The monodromy group $\Mon_I$ is the group generated by
the monodromy maps associated with the Gauss-Manin
connections for all families of deformations of $(M,I)$.

\hfill

\remark
The intersection of the monodromy group and
the Torelli group is finite (\cite[Corollary 7.3]{_V:Torelli_}).

\hfill

Using the erroneous quotation from
\cite{_Huybrechts:finiteness_},
I claimed that the Teichm\"uller space 
has finitely many connected components
(\cite[Corollary 1.12]{_V:Torelli_}.
Generally speaking, this is false.
Indeed, as follows from
\ref{_Torelli_on_components_Theorem_} below,
the Torelli group acts on the connected components
of the Teichm\"uller space with finite stabilizers
and finitely many orbits. Therefore, the 
Teichm\"uller space has infinitely many 
connected components if and only if
the Torelli group is infinite.

\cite[Corollary 1.12]{_V:Torelli_} was used to prove the
the following theorem.

\hfill

\theorem\label{_mono_finite_in_O(H^2)_Theorem_}
Let $\Mon_I$ be the monodromy group of a hyperk\"ahler
manifold $M$, and $\Phi:\; \Gamma\arrow O(H^2(M,\Z))$ the
natural homomorphism. Then 
\begin{description}
\item[(i)] The group $\Phi(\Mon_I)$ has
finite index in $O(H^2(M,\Z))$, for all $I$.
\item[(ii)] The number of connected components of 
the marked moduli space of $M$ is finite.
\end{description}

\proof
The proof is similar the one given by
Bakker and Lehn (\cite[Theorem 8.2]{_Bakker_Lehn:singular_}),
where they prove this result for projective deformations
of singular hyperk\"ahler manifolds. For non-algebraic
deformations, their proof needs some adjustments.


First, notice that
(i) is equivalent to (ii). Indeed, consider the action of the group
$\Gamma/K$ on ${\cal M}=\Teich/K$. The space ${\cal M}$ has finitely many
components if and only if the stabilizer $\tilde \Mon_I$ of each connected component of 
${\cal M}$ in $\Gamma$ has finite
index in $\Gamma/K$. However, the group $\tilde \Mon_I$ is generated by
$K$ and $\Mon_I$, while $\Phi(K)=1$ by definition of $K$. Therefore
$\Phi(\tilde \Mon_I)=\Phi(\Mon_I)\subset \Phi(\Gamma)$. The group $\Phi(\Gamma)$
has finite index in $O(H^2(M,\Z))$ by 
\ref{_mapping_class_intro_Equation_}.

Now we prove \ref{_mono_finite_in_O(H^2)_Theorem_} (i).
Fix a very ample line bundle $L$ on $(M,I)$,
and let $\tilde Z\subset {\cal M}$ be the set of all 
complex structures $J$ on $M$
such that $c_1(L)$ is very ample on $(M,J)$.
All such manifolds $(M,J)$ can be embedded into ${\Bbb P}^N$
with a fixed $N$ which is determined from
the Riemann-Roch formula. Since the Hilbert
scheme with prescribed Poincar\'e
polynomial is bounded, $Z$ is a union of 
finitely many deformation families (\cite{_Matsusaka:polarized_}).

Consider the Teichm\"uller space $\Teich_L$ of all
$J\in \Teich$ such that $\eta=c_1(L)\in H^2(M,\Z)$ is the first
Chern class of a very ample line bundle on $(M,J)$.
Then $\eta \in H^{1,1}_J(M)$, and $H^{1,1}_J(M)$
can be obtained as the orthogonal complement of
the 2-plane $l_J:=\langle \Re\Omega_J, \Im\Omega_J\rangle$,
where $\Omega_J$ is cohomology class of a 
holomorphically symplectic form on $(M,J)$.
The corresponding period space $\Perspace_\eta$ is
the Grassmannian $\Gr_{+,+}(\eta^\bot)$.
Since $\eta^\bot$ has signature $(2, b_2-3)$,
the space
\[ 
\Perspace_\eta=\Gr_{+,+}(\eta^\bot)=\frac{SO^+(2, b_2-3)}{SO(2)\times SO(b_2-3)}
\]
is a symmetric space, with the arithmetic group
$SO(\eta^\bot \cap H^2(M,\Z))$ acting   properly by isometries.
Let $\Gamma_\eta$ be a subgroup of all elements of $\Gamma$ fixing $\eta$.
Since $\Phi(\Gamma)$ is a finite index subgroup in $O(H^2(M,\Z))$, 
$\Phi(\Gamma_\eta)$ is a finite index subgroup in $O(\eta^\bot \cap H^2(M,\Z))$.
Therefore, the quotient space $\Perspace_\eta/\Gamma_\eta$ is quasiprojective
by the Baily-Borel theorem (\cite{_Baily_Borel:ann_}).

Let $\Teich_\eta$ be the Teichm\"uller space of all $J$ such that
$\eta$ is of type (1,1) on $(M,J)$. 
For any point $J \in \Teich_\eta$ with 
the Picard group of $(M,J)$ has rank 1, the bundle
$L$ is ample (\cite{_Huybrechts:cone_}). Replacing $L$ by $L^N$ if necessary,
where $N$ is the appropriate Fujita constant, we
may assume that $L$ is very ample. 
This implies that $\Teich_L$ is dense 
in $\Teich_\eta$ and has the same number
of connected components.

Consider the Hilbert scheme $Z$ defined above. Its universal
cover is mapped to $\Perspace_\eta$, because $\Perspace_\eta$
is the classifying space of the Hodge structures on $H^2(M)$.
This gives a complex analytic map $Z\arrow \Perspace_\eta/\Phi(\Gamma_\eta)$,
which is algebraic by Borel's extension theorem
(\cite[Theorem 3.10]{_Borel:extension_}). Therefore, a general point
in $\Perspace_\eta/\Phi(\Gamma_\eta)$ has finite preimage in
$Z= \Teich_L/\Gamma_\eta$. On the other hand,
each connected component of $Z$ can be identified
with a connected component of $\Teich_L/\Mon_\eta^I$,
where $\Mon_\eta^I$ is the corresponding subgroup of
the monodromy group. We obtain that each component of $\Teich_L/\Mon_\eta^I$
maps to $\Perspace_\eta/\Phi(\Gamma_\eta)$ with finite fibers,
and $\Phi(\Mon_\eta^I)$ has finite index in $\Phi(\Gamma_\eta)$.

The union
$\bigcup_\eta \Phi(\Mon_\eta^I)$ generates a finite
index subgroup in $O(H^2(M,\Z)))$,
as follows  from \ref{_codim1_f_index_Lemma_} below.
 Therefore, $\Phi(\Mon_I)$ also has
finite index in $O(H^2(M,\Z))$. 
\endproof

\hfill 

\lemma\label{_codim1_f_index_Lemma_}
Let $(\Lambda, q)$ be a non-degenerate quadratic lattice of signature 
$(p, q)$, $p>2$ and $q>1$, 
$O(\Lambda)$ its isometry group, and $O_\eta(\Lambda)$
its subgroup fixing
$\eta\in \Lambda$ which satisfies $c:=(\eta,\eta)>0$.
For each $\eta\in \Lambda$ with positive square, fix
a subgroup $\Gamma_\eta\subset O_\eta(\Lambda)$ of finite index.
Let $S:=\bigcup_\eta\Gamma_\eta \subset O(\Lambda)$
be the union of all $\Gamma_\eta$ for all such $\eta$.
Then $S$ generates a finite index subgroup $\Gamma_S\subset O(\Lambda)$.

\hfill

{\bf Proof. Step 1:}
Let $V= \R^{p,q}:=\Lambda\otimes_\Z \R$
be the vector space associated with $\Lambda$,
and  $X_c\subset V$ the set of all vectors of square $c=(\eta, \eta)$.
We prove that $\Gamma_S$ acts on the quadric $X_c$ ergodically
for a (unique) $O(V)$-invariant Lebesgue measure.

Let $Y_r$ be the set of 
all $x\in X_c$ such that $g(x, \eta)=r$.
The space $Y_r$ is equipped with a homogeneous
action of the group $H:= O(\eta^\bot\otimes_\Z \R)=O(p-1, q)$,
and the stabilizer $H_0$ of a generic point $z\in X_c$ is either $O(p-2, q)$ or
$O(p-1, q-1)$, depending on the signature of the
space $\langle x, \eta\rangle$. Therefore, $H_0$ is non-compact.

Applying Moore's ergodicity theorem (\cite[Theorem 7]{_Moore:ergodi_})
to the space $Y_r=H/H_0$, we see that
$\Gamma_\eta$ acts ergodically on  $Y_r$. Therefore, any
$\Gamma_\eta$-invariant measurable function on $X_c$
is $H$-invariant almost
everywhere. Applying this to two non-collinear vectors
$\eta_1, \eta_2$ with positive square, we obtain that
any $\Gamma_S$-invariant measurable function on $X_c$ is
both $O(\eta_1^\bot\otimes_\Z \R)$-invariant and
$O(\eta_2^\bot\otimes_\Z \R)$-invariant, hence constant.

\hfill

{\bf Step 2:} We are going to prove that 
the group $\Gamma_S$ has finite covolume in $O(V)$.
Consider the map $O(V) \arrow X_c$ taking $g\in O(V)$ to
$g(\eta)$, where $\eta\in \Lambda$.
This is a locally trivial fibration
with fiber $H= O(\eta^\bot\otimes_\Z \R)$. Choose an open set of finite measure 
$U\subset X_c$ such that $\Gamma_S U\supset X_c$
(Step 1).  Since $\Gamma_\eta\subset H$
 is an arithmetic lattice in $H$ (\cite{_Morris_arithm_}), 
there exists a subset  $V\subset H$ of finite measure such that 
$\Gamma_\eta V\supset H$.
Choose a section $R\subset O(V)$ of the map $O(V) \arrow X_c$ 
over $U\subset X_c$, and let $W\subset O(V)$ be obtained
as $V\cdot R$. The set $W= U\times V$ has finite measure
by Fubini's theorem. To prove that $\Gamma_S$ has finite
covolume in $O(V)$, it remains to show that 
$\Gamma_S W\supset O(V)$. Since $\Gamma_\eta V\supset H$,
we have $\Gamma_S W\supset\pi^{-1}(U)$. Since $\Gamma_S \supset X_c$,
we also have 
\[ \Gamma_S W\supset\Gamma_S\pi^{-1}(U)=\pi^{-1}(\Gamma_SU)=O(V).
\]
Then $\Gamma_S$ has finite covolume in $O(V)$.

\hfill

{\bf Step 3:} The number $\frac{\Vol(O(V)/\Gamma_S)}{\Vol(O(V)/O(\Lambda))}$
is equal to the index of $\Gamma_S$ in $O(\Lambda)$; Step 2 implies
that this number is finite.
\endproof

\section{The Torelli group action on the Teichm\"uller space}
\label{_MCG_on_Teich_Section_}

In \cite[Remark 1.12]{_V:Torelli_} I 
misquoted a paper \cite{_Huybrechts:finiteness_} of Dan Huybrechts 
stating that the Teichm\"uller space of a hyperk\"ahler
manifold has finitely many components. This is false,
as explained above. The correct version
of this statement is as follows.

\hfill

\theorem\label{_Torelli_on_components_Theorem_}
Let $M$ be a compact hyperk\"ahler manifold,
$\Gamma$ its mapping class group, $\Teich$ the
Teichm\"uller space of complex structures of hyperk\"ahler
type, and $K$ the Torelli group, that is, the group
of all elements $\gamma\in \Gamma$ acting trivially on
$H^2(M, \R)$. Then 
\begin{description}
\item[(i)]
$K$ acts on the space of connected
components of $\Teich$ with finitely many
orbits, and the stabilizer of each component is finite. 
\item[(ii)]
Moreover, for each $\gamma \in K$ fixing a point
$I\in \Teich$, $\gamma$ acts trivially on the
connected component $\Teich^I$ of $I$ in $\Teich$.
\end{description}

\proof
\ref{_Torelli_on_components_Theorem_} (i) is 
\ref{_mono_finite_in_O(H^2)_Theorem_} (ii) and
\ref{_Torelli_on_components_Theorem_} (ii) is 
\ref{_Torelli_action_Claim_}.
\endproof

\hfill

\corollary The Torelli space of complex structures
of hyperk\"ahler type on $M$ has infinitely many connected 
components if the Torelli group is infinite, and finitely many
components if it is finite.
\endproof

\hfill

\remark
A similar result is true for the Teichm\"uller space
of complex structures of K\"ahler type on 
a compact torus, which also has infinite Torelli group
(\cite{_Hatcher:torus_}). The Albanese map from a compact
K\"ahler torus
to its Albanese space is an isomorphism, and this
provides the torus with the canonical flat coordinates.
Then the connected component of its Teichmuller
space is the set $GL^+(2n, \R)/GL(n, \C)$ of complex
structures on an oriented vector space. The Torelli group
acts on the set of connected components freely and 
transitively for the same reason as above. Indeed,  
an element $\gamma$ of the Torelli group which stabilizes
$I\in \Teich$ commutes with the Albanese map, and
hence acts trivially on the whole connected 
component $\Teich_I$ of the Teichm\"uller space. Then 
$\gamma$ defines an automorphism of $(M,J)$ for
all $J\in \Teich_I$. For a generic complex torus, all automorphisms
are homotopic to identity, hence $\gamma$ is trivial.

\section{Errata}

The following statements of \cite{_V:Torelli_}
are false and should be amended: Remark 1.12, Theorem 1.16,
Theorem 3.4, Theorem 3.5 (iv), Theorem 4.26 (ii) and (iii).
Correct versions of these results are given 
in Section \ref{_MCG_Section_} (Theorem 3.4, Theorem 3.5. Theorem 1.16),
Section \ref{_marked_Section_} (Remark 1.12), Section \ref{_MCG_on_Teich_Section_}
(Theorem 4.26 (ii) and (iii)). 
These erroneous claims were used in the proof
of \cite[Theorem 4.25, Corollary 4.31, Corollary 7.3]{_V:Torelli_}.
\ref{_Torelli_on_components_Theorem_} is sufficient
to prove these results.

\hfill

\remark
Using this opportunity to update the paper,
I want to address a minor gap in the proof of
Corollary 7.3, which claims ``Using \cite[Remark 13]{_Catanese:moduli_}, we may
assume that there exists a universal fibration
$\pi:Z\arrow \Teich_I$.'' In a more recent
paper \cite{_Markman:universal_}, 
E. Markman gives an explicit construction of the 
universal fibration.

\hfill

Summing it up: due to my oversight, 
some claims about the Teichm\"uller spaces made in
\cite{_V:Torelli_} were incorrect, and some results 
of D. Huybrechts and D. Sullivan were misquoted.
Most of the errors can be corrected if we replace the 
Teichm\"uller space by the marked moduli space, and
the mapping class group $\Gamma$ by its quotient
$\Gamma/K_0$, where $K_0$ is a subgroup 
of all elements acting trivially on ${\cal M}$,
which has finite index in the Torelli group.

\hfill

{\bf Acknowledgements:}
I am grateful to Dennis Sullivan, Andrey Soldatenkov, 
Matthias Kreck and Richard Hain for calling my attention
on the erroneous claims in \cite{_V:Torelli_}. Apologies to those who were
mislead by my errors. Many thanks to Dmitry Kaledin and
Benjamin Bakker for helpful remarks and finding inaccuracies 
in the earlier versions of this erratum.

\hfill

{\scriptsize

\noindent
{\sc Misha Verbitsky\\
 Instituto Nacional de Matem\'atica Pura e
              Aplicada (IMPA) \\ Estrada Dona Castorina, 110\\
Jardim Bot\^anico, CEP 22460-320\\
Rio de Janeiro, RJ - Brasil \\
also:\\
\sc Laboratory of Algebraic Geometry,\\
National Research University HSE,\\
Department of Mathematics, 6 Usacheva Str. Moscow, Russia\\
\tt  verbit@impa.br}.
}

\end{document}